%
%
%
%
\input amstex
\documentstyle{amsppt}
\vcorrection{-7mm} \hcorrection{-10mm}
 \NoBlackBoxes \raggedbottom
\pageheight{9.0 in} \pagewidth{150mm}
 \topmatter
\smallskip

\NoBlackBoxes
\topmatter
\title   Browder-Livesay filtrations and  the example of Cappell and Shaneson
\endtitle
\rightheadtext{On the example of Cappell and Shaneson}
\author
F. Hegenbarth, 
Y. V. Muranov, and  D. Repov\v s
\endauthor
\thanks {\it Acknowledgements.} The second author was supported by the  CONACyT Grant 98697 and he  thanks the Commission on 
Development and Exchange of the International Mathematical Union for the 
travel support. The third author was supported by the SRA grants P1-0292-0101, J1-2057-0101, and J1-4144-0101.
 \endthanks
\abstract 
Let $M^3$ be a 3-dimensional manifold with fundamental group 
$\pi_1(M)$ which  contains a  quaternion subgroup $Q$ of order 8.
In 1979 Cappell and Shaneson constructed a nontrivial normal map  
$
f\colon M^3\times T^2\to M^3\times S^2$ 
which  cannot be  detected  by simply connected surgery obstructions along submanifolds  of codimension 0, 1, or 2, but it can be  detected  by the codimension 3 Kervaire-Arf invariant.  
   The proof of non-triviality of  $\sigma(f)\in L_5(\pi_1(M))$ is based 
on consideration of a Browder-Livesay filtration  of a manifold $X$ with 
$\pi_1(X)\cong \pi_1(M)$.  
For a Browder-Livesay pair $Y^{n-1}\subset X^n$, the restriction of a normal  map to the submanifold $Y$ is given by a partial multivalued map 
$\Gamma\colon L_n(\pi_1(X))\to L_{n-1}(\pi_1(Y))$, and the  Browder-Livesay filtration provides an  iteration $\Gamma^n$. This map is a basic step in the definition of the iterated Browder-Livesay invariants 
which give obstructions to realization of surgery obstructions by normal maps of closed manifolds. 

In the present paper we prove that $\Gamma^3(\sigma(f))=0$  for any Browder-Livesay filtration of a manifold $X^{4k+1}$ with $\pi_1(X)\cong Q$.
 We compute   splitting obstruction groups 
for various inclusions  $\rho\to Q$ of index 2, describe  natural maps in the braids of exact sequences, and  make more precise several results about  surgery obstruction groups of the group $Q$. 
\endabstract

\endtopmatter

\document

\subhead 1. Introduction 
\endsubhead
Let $M^3$ be a 3-dimensional manifold with fundamental group 
$\pi_1(M)$ which  contains the   quaternion subgroup 
$$
Q=\{x,y|x^4=y^2, yxy^{-1}=x^{-1}\}
\tag 1.1
$$
of order 8. Consider  a map 
$$
f\colon M^3\times T^2\longrightarrow M^3\times S^2
\tag 1.2
$$
where $T^2=S^1\times S^1$ is the torus equipped with the Lie group invariant framing \cite{1}, and 
$
T^2 \longrightarrow S^2
$ 
is a normal  map with a nontrivial Kervaire-Arf invariant. 
In 1979 Cappell and Shaneson \cite{1} proved that surgery obstruction 
$\sigma(f)$ is nontrivial.
The nontriviality of $\sigma(f)$ cannot be detected by simply connected surgery obstructions along submanifolds of codimension 0, 1, and 2, but it can be detected by the codimension 3 Arf-invariant (see \cite{1} and \cite{5}).

 Recall that a pair of manifolds 
$$
Y^{n-1}\subset X^n
\tag 1.3 
$$ 
is called {\it a Browder-Livesay pair}  if $Y$ is a codimension one locally-flat closed submanifold and the natural inclusion 
$\pi_1(Y)\to \pi_1(X)$ induces an isomorphism of fundamental groups. 
In what follows we shall consider only the case when 
$\operatorname{dim} Y=n-1\geq 5$.  
For a Browder-Livesy pair, we can write down a push-out  square of fundamental groups 
$$
\CD 
\pi_1(\partial U)@>\cong >> \pi_1(X\setminus Y)\\
@VVV @VVV \\
\pi_1(Y)@>>> \pi_1(X)\\
\endCD\ 
= \
\CD 
\rho@>>> \rho\\
@VVi_- V @VVi_+ V \\
G^-@>>> G^+\\
\endCD\tag 1.4
$$
where  $\partial U$ is the boundary of a tubular neighborhood of $Y$ in $X$,  and $X\setminus Y$ is the closure of the complement of a tubular neighborhood. The groups in square (1.4) are equipped with orientation. The upper horizontal map  and 
 the vertical maps  agree with orientations. The bottom horizontal
 map preserves the orientation on the images of the vertical maps which are inclusions of index 2  and
 reverses orientation outside of these images.   We shall denote
this fact by superscript "$+$" or "$-$". We shall omit this
superscript if the orientation is clear from the context.

For a Browder-Livesay pair (1.3) we have a braid of exact sequences (see \cite{3} and \cite{8}):
$$
 \matrix \rightarrow & {L}_{n}(\rho) & \overset{i_*}\to{\longrightarrow} &
   L_{n}(G) &
\overset{\partial}\to{\rightarrow} & LN_{n-2}(\rho\to
G)&\rightarrow
 \cr
\ & \nearrow \  \  \ \ \ \ \  \ \searrow &\ &s\nearrow \ \ \ \ \
                           \searrow
   & \ & \nearrow \ \ \  \ \ \ \ \ \ \searrow & \ \cr
  \ & \ & LP_{n-1}(F)& \Gamma\downarrow & L_{n}(i_*) & \ & \ \cr
   \ & \searrow \ \ \ \ \ \ \ \ \ \nearrow &\ &q\searrow  \ \ \ \ \
\nearrow
 & \ & \searrow \ \   \ \ \ \  \ \ \ \nearrow & \ \cr
\rightarrow & LN_{n-1}(\rho\to G) &
\overset{c}\to{\longrightarrow} &
 L_{n-1}(G^-) &
 \overset{i^!_-}\to{\longrightarrow} & {L}_{n-1}(\rho) & \rightarrow
\endmatrix
\tag 1.5
 $$
where $LP_{n-1}(F)\cong L_n(i^!_-)$ are  surgery
obstruction groups for the manifold pair $(X,Y)$,  $LN_{n}(\rho\to G)$ are the splitting obstruction groups,  and  $L_n(i_*)$ are the relative surgery
obstruction groups for the inclusion map $i$ (see \cite{3}, \cite{8} and 
\cite{9}). The horizontal rows in (1.5)   are chain complexes and  $\Gamma$ is an  isomorphism of the
corresponding homology groups. The maps $s$ and $q$ are
the natural forgetful maps,  and  the map $c$ denotes  passing from
surgery problem inside the manifold $X^n$ to an abstract surgery
problem. For a normal map $f\colon M\to X$, the restriction to a transversal preimage 
$$
f|_{f^{-1}(Y)}\colon f^{-1}(Y)\to Y^{n-1}
$$
is a normal map to $Y^{n-1}$, and  the corresponding map of surgery obstruction groups  is given by a partial multivalued map 
$\Gamma\colon L_n(\pi_1(X))\to L_{n-1}(\pi_1(Y))$.

Let  $\Cal X$ be a filtration
$$
X_k\subset X_{k-1}\subset \cdots \subset X_2\subset X_1\subset X_0=X
\tag 1.6
$$
of a closed $n$-dimensional manifold $X$ by means of locally flat closed
submanifolds.  A filtration in (1.6) is called a {\it  Browder-Livesay filtration}, if 
every pair of submanifolds $(X_i,X_{i+1})\ (0\leq i\leq k-1) $ is a
Browder-Livesay pair 
(see \cite{2}, \cite{5} and \cite{7}). In what follows we shall  assume that
$\text{dim}\, \,  X_k=n-k\geq 5$. 

Let  $F_i \ (0\leq i\leq k-1)$ be a square of fundamental groups in
the splitting problem for a manifold pair $(X_i,X_{i+1})$  of a
filtration in  (1.6),  $G_i=\pi_1(X_i)$, and
$\rho_i=\pi_1(\partial U_{i+1})\cong \pi_1(X_i\setminus X_{i+1})$, 
where $\partial U_{i+1}$ is a boundary of a tubular neighborhood of
$X_{i+1}$ in $X_{i}$. 
 Then  $LN_*(\rho_i\to G_i)$
are the splitting obstruction groups for the manifold pair
$(X_i,X_{i+1})$.

 Every inclusion $\rho_i\to G_i$ of index 2 of oriented groups gives a commutative braid of exact
 sequences   (1.5).
Putting together central squares from these diagrams (see \cite{4}, \cite{5} and  \cite{7} ) we obtain a
  commutative diagram

$$
\matrix
 \overset{}\to{\longrightarrow} &   L_{n}(G_0) & \overset{\partial_0}\to{\rightarrow} &LN_{n-2}(\rho_0\to G_0) \cr
   &s\nearrow \ \ \ \ \                           \searrow p &  &\cr
  LP_{n-1}(F_0)& \Gamma\downarrow & L_{n}(\rho_0\to G_0) &\cr
    &q\searrow \ \ \ \ \   \nearrow  r &  &\cr
\longrightarrow &   L_{n-1}(G_1) &
\overset{\partial_1}\to{\rightarrow}  &LN_{n-3}(\rho_1\to
G_1)\cr
   &s\nearrow \ \ \ \ \                           \searrow p &  &\cr
  LP_{n-2}(F_1)& \Gamma\downarrow & L_{n-1}(\rho_1\to G_1) &\cr
    &q\searrow \ \ \ \ \   \nearrow r &  &\cr
\longrightarrow &   L_{n-2}(G_2) &
\overset{\partial_2}\to{\rightarrow} &LN_{n}(\rho_2\to G_2) \cr
   &s\nearrow \ \ \ \ \                           \searrow p & & \cr
  LP_{n-3}(F_2)& \Gamma\downarrow & L_{n-2}(\rho_2\to G_2) &\cr
    &q\searrow \ \ \ \ \   \nearrow r & & \cr
& L_{n-3}(G_3) & &\cr
                           & \vdots&&& &&\cr
\longrightarrow &   L_{n-k+1}(G_{k-1}) &
\overset{\partial_{k-1}}\to{\rightarrow} &LN_{n-k-1}(\rho_{k-1}\to G_{k-1}) \cr
   &s\nearrow \ \ \ \ \                           \searrow p & & \cr
  LP_{n-k}(F_{k-1})& \Gamma\downarrow & L_{n-k+1}(\rho_{k-1}\to G_{k-1}) &\cr
    &q\searrow \ \ \ \ \   \nearrow r & & \cr
& L_{n-k}(G_k). & &\cr
 \endmatrix
\tag 1.7
 $$
 \smallskip

\noindent In this diagram we denote by $s$, $q$, $p$,  and $r$   similar maps from
different diagrams. However, in what follows,  it will be clear from
the context which map is under consideration. Note that the groups
and the maps in diagram (1.7) are defined by the subscripts taken
$\bmod \ 4$.

Now we recall  an inductive definition of  the sets
$$
\Gamma^j(x)\subset L_{n-j}(G_j) \ \text{for} \ (0\leq j\leq k)
$$
and the iterated Browder-Livesay $j$-invariants $(1\leq j\leq k)$ with
respect to  filtration  (1.6) (see \cite{5} and \cite{7}).

\proclaim{Definition 1}  Let  $x\in L_n(G_0)$. By definition,
$$
\Gamma^0(x)=\{x\}\subset L_n(G_0).
$$
The set $\Gamma^0(x)$ is said to be  trivial if $x\in
\operatorname{Image}\{L_n(\rho_0)\to L_n(G_0)\}$.
 Let a set
 $$
 \Gamma^{j}(x)\subset L_{n-j}(G_{j})\ \ (0\leq j\leq k-1)
 $$
 be defined. For $j\geq 1$, it is called trivial if $0\in \Gamma^{j}(x)$.

If $\Gamma^{j}(x) (0\leq j\leq k-1)$ is defined and nontrivial, then the $(j+1)$-th
Browder-Livesay invariant with respect filtration  (1.6) is the set
$$
\partial_{j} (\Gamma^{j}(x))\subset LN_{n-j-2}(\rho_{j-1}\to
G_{j-1}).
$$
 The
$(j+1)$-th invariant is nontrivial if $0\notin \partial_{j}
(\Gamma^{j}(x))$.

If the $(j+1)-th \ (1\leq j\leq k-1)$ Browder-Livesay invariant is
defined and trivial then the set  $\Gamma^{j+1}(x)$ is defined as
$$
\Gamma^{j+1}(x)\overset{def}\to{=}
\Gamma(\Gamma^{j}(x))\overset{def}\to{=} \{qs^{-1}(z)|z\in
\Gamma^{j}(x), \partial_{j}(z)=0\}\subset L_{n-j-1}(G_{j+1}).
$$
\endproclaim

\proclaim{Theorem 1} (See \cite{5}.)  Let $x\in L_n(G_0)$
be an element with  a  nontrivial $j$-th Browder-Livesay  invariant
for some $j\geq 1$ relative to a Browder-Livesay filtration $\Cal
X$ of the manifold $X$. Then the element $x$ cannot be realized by a
normal map of closed manifolds.
\endproclaim

Note that the necessary condition for nontriviality of the  $j$-th Browder-Livesay  invariant  of $x\in L_n(G_0)$ $(1\leq j)$ is 
nontriviality of $\Gamma^{j-1}(x)$. 

Taking in (1.2) a 3-dimensional manifold $M$ with $\pi_1(M)\cong Q$ we obtain 
a surgery obstruction $\sigma(f)\in L_1(Q)$ (see \cite{1}).   The proof of nontriviality  $\sigma(f)$  in \cite{1} is based 
on a consideration  of the Browder-Livesay filtration   
$$
X_3\subset X_2\subset X_1\subset X_0=X, \ \ (\operatorname{dim} X_3\geq 5)
\tag 1.8
$$
with $\pi_1(X)=Q$, $\pi_1(X_1)=Q^{+,-}$, $\pi_1(X_2)=Q^{-,-}$, and 
$\pi_1(X_3)=Q$. 
In particular, it follows from \cite{1}, that  $\Gamma^3(\sigma(f))=0$
for  filtration in (1.8).   

 In  \cite{5} was  introduced the notion of a type of any element  $x\in L_n(\pi)$ and it was proved  that the elements of the first and the second type cannot be realized by normal maps of closed manifolds. For any  element $x$ of  the second type there exists an "infinite" Browder-Livesay filtration for  which $\Gamma^k(x)$ is nontrivial for all $k\geq 0$. Every element $x\in L_n(\pi)$ which lays in the subgroup generated by surgery obstructions of normal maps of closed manifolds has the third type. And  if an element $x$ has  the third type, then for any Browder-Livesay filtration there exists 
 a finite $k$, such that $\Gamma^k(x)$ is trivial. In all known to the  authors cases $\Gamma^3(x)$ is trivial for the elements of the third type. 
 
 In the present paper,   we prove that $\Gamma^3(\sigma(f))$ is trivial for any Browder-Livesay filtration of a manifold $X^{4k+1}$ with $\pi_1(X)\cong Q$,  where  $\sigma(f)\in L_1(Q)$. 
 We compute   splitting obstruction groups $LN_*$ and surgery obstruction groups $LP_*$ 
for various inclusions  $\rho\to Q$ of index 2, describe  natural maps in the braids of exact sequences in (1.5), and  make more precise several results about  surgery obstruction  groups of the group $Q$. 
 
We use the  surgery and splitting obstruction groups equipped with decoration "s" and we do not mention this in designations. 

\smallskip

\subhead 2. Browder-Livesay filtrations for a manifold $X^{4k+1}$ with 
$\pi_1(X)=Q$
\endsubhead
\bigskip

Let $Q$ be the quaternion group from (1.1).
We define an  orientation homomorphism 
$w: Q\to \{\pm 1\}$   on generators $x$ and $y$ and  denote it by  superscripts  in the following way:
$$
w(x)=w(y)=1, \ \ (Q,w)=Q^+;
$$
$$
 w(x)=1, w(y)=-1, \ \ (Q,w)=Q^{+,-};
$$
$$
w(x)=-1, w(y)=1, \ \ (Q,w)=Q^{-,+};
$$
and
$$ 
w(x)=-1, w(y)=-1, \ \ (Q,w)=Q^{-,-}\ \ \ \text{(in this case $w(xy)=1$)}.
$$

Let $\rho=\Bbb Z/4$ be the cyclic group with generator $t$. There are only two   homomorphisms of orientation  on the group $\rho$ (homomorphisms $\Bbb Z/4$ into the group ${\pm 1}$), and we shall write $\rho^+=\Bbb Z/4^+$ in the case of trivial orientation and 
$\rho^-=\Bbb Z/4^-$ in the opposite case. The group $Q$ has only three different subgroups of index 2 generated by $x$, $y$, and $xy$. All these 
subgroups are isomorphic to $\Bbb Z/4$.

Consider the Browder-Livesay filtration 
$$
X_3\subset X_2\subset X_1\subset X_0=X, \ \ (\operatorname{dim} X_3\geq 5)
\tag 2.1
$$
of a manifold $X$ with  $\pi_1(X)
\cong Q$.

Filtration in (2.1) yields
 a commutative diagram  
$$
\matrix
        &         &         & \rho_2&         &   &         & \rho_1&
              &     &         &\rho_0&         &   \\
    &         &\swarrow &       &\searrow & &\swarrow &
     &\searrow &     &\swarrow &      &\searrow &     \\
        & \pi_1(X_3)&         &  \overset{\cong}\to{\to}     &         &\pi_1(X_2)&
               &  \overset{\cong}\to{\to}      &         & \pi_1(X_1) &
                      & \overset{\cong}\to{\to}     &         &\pi_1(X_0)\\
\endmatrix
\tag 2.2
$$
where  every triangle 
corresponds to the square $F_i$, since  we have  isomorphisms of oriented  groups $\rho_i=\pi_1(\partial U_{i+1})\cong \pi_1(X_{i}\setminus X_{i+1})$ for $0\leq i\leq 2$. 

 The diagram in (2.2) is  commutative as a diagram of groups,
and the skew  maps are inclusions of index 2 preserving the
orientations. 
 The horizontal maps in (2.2) 
preserve the orientation on the images of skew  maps  and reverse the
orientation outside these images. 
Two  diagrams as in (2.2) are isomorphic if there exists an isomorphism between them preserving orientations.

\proclaim{Theorem 2} Let $X^{n}$, $n\geq 8$,  be a closed topological  manifold with $\pi_1(X)\cong Q$.  Then the diagram (2.2) is isomorphic to   one of the following diagrams  
$$
\matrix
        &         &         & \Bbb Z/4^+&         &   &         & \Bbb Z/4^-&
              &     &         &\Bbb Z/4^+&         &   \\
    &         &^{k}\swarrow &       &\searrow^{k} & &^{j}\swarrow &
     &\searrow^{j} &     &^{i}\swarrow &      &\searrow ^{i}&     \\
        & Q^{-,+}&         &  \overset{=}\to{\longrightarrow}     &         &Q^{-,-}&
               &  \overset{=}\to{\longrightarrow}      &         & Q^{+,-} &
                      & \overset{=}\to{\longrightarrow}     &         &Q^{+}\\
\endmatrix
\tag 2.3
$$
where skew homomorphisms  are given on the generator $t$  by the maps $i(t)=x, \ j(t)=y,\  k(t)=xy$; 
$$
\matrix
        &         &         & \Bbb Z/4^+&         &   &         & \Bbb Z/4^+&
              &     &         &\Bbb Z/4^+&         &   \\
    &         &^{k}\swarrow &       &\searrow^{k} & &\swarrow &
     &\searrow &     &\swarrow &      &\searrow &     \\
        & Q^{-,+}&         &  \overset{=}\to{\longrightarrow}     &         &Q^{+}&
               &  \overset{=}\to{\longrightarrow}      &         & Q^{+,-} &
                      & \overset{=}\to{\longrightarrow}     &         &Q^{+}\\
\endmatrix
\tag 2.4
$$
where  $k:\Bbb Z/4\to Q$ is given on the generator $t$ by the map $t\to y$ and other skew  homomorphisms  are given on the generator $t$  by the map $t\to x$; 

$$
\matrix
        &         &         & \Bbb Z/4^+&         &   &         & \Bbb Z/4^+&
              &     &         &\Bbb Z/4^+&         &   \\
    &         &\swarrow &       &\searrow & &\swarrow &
     &\searrow &     &\swarrow &      &\searrow &     \\
        & Q^{+,-}&         &  \overset{=}\to{\longrightarrow}     &         &Q^{+}&
               &  \overset{=}\to{\longrightarrow}      &         & Q^{+,-} &
                      & \overset{=}\to{\longrightarrow}     &         &Q^{+}\\
\endmatrix
\tag 2.5
$$
where all homomorphisms   $\Bbb Z/4\to Q$ are given on the  generator $t$ by the map $t\to x$;

$$
\matrix
        &         &         & \Bbb Z/4^-&         &   &         & \Bbb Z/4^-&
              &     &         &\Bbb Z/4^+&         &   \\
    &         &^{k}\swarrow &       &\searrow^{k} & &^{j}\swarrow &
     &\searrow^{j} &     &^{i}\swarrow &      &\searrow^{i}&     \\
        & Q^{-,+}&         &  \overset{=}\to{\longrightarrow}     &         &Q^{-,-}&
               &  \overset{=}\to{\longrightarrow}      &         & Q^{+,-} &
                      & \overset{=}\to{\longrightarrow}     &         &Q^{+}\\
\endmatrix
\tag 2.6
$$
where  skew  homomorphisms  are given on the generator $t$  by the maps $i(t)=x, \ j(t)=y,\  k(t)=y$; 

$$
\matrix
        &         &         & \Bbb Z/4^-&         &   &         & \Bbb Z/4^-&
              &     &         &\Bbb Z/4^+&         &   \\
    &         &^{k}\swarrow &       &\searrow^{k} & &^{j}\swarrow &
     &\searrow^{j} &     &^{i}\swarrow &      &\searrow ^{i}&     \\
        & Q^{-,+}&         &  \overset{=}\to{\longrightarrow}     &         &Q^{-,-}&
               &  \overset{=}\to{\longrightarrow}      &         & Q^{+,-} &
                      & \overset{=}\to{\longrightarrow}     &         &Q^{+}\\
\endmatrix
\tag 2.7
$$
where  skew  homomorphisms  are given on the generator $t$  by the maps $i(t)=x, \ j(t)=y,\  k(t)=x$. 

For any diagram from the list (2.3)--(2.7) there exists a Browder-Livesay filtration (2.1) of the manifold $X^n$ which gives this diagram. \endproclaim
\demo{Proof} There exist isomorphisms of the oriented  groups
$$
Q^{+,-}\cong Q^{-,+}\cong Q^{-,-}
$$
since there exist automorphisms of the group $Q$ permuting $x$, $y$,
and $xy$.  We have only three index 2 subgroups  of $Q$
generated by $x$, $y$, and $xy$, which   are isomorphic to $\Bbb Z/4$. From this follows that all skew maps in (2.3) are inclusions  of
$\Bbb Z/4$ to $Q$ preserving orientation. 
Now a consideration of various cases which agree with orientation conditions on the horizontal maps gives the first statement of the proposition. Consider 
the characteristic map $\phi\colon X^n\to RP^N$ ($N$ sufficiently large) 
of the subgroup $i(\Bbb Z/4^+)=i(\rho^+)\subset Q^+$ generated by $x$  such that 
$\operatorname{Ker} \phi =\rho^+$. The transversal preimage of $RP^{N-1}$
gives a submanifold  $X_1\subset X$ such that $X_1\subset X$ is a Browder-Livesay pair with the commutative triangle 
$$
\matrix 
 &        & \rho^+ & & \\
&\swarrow &      &\searrow& \\
Q^+\cong \pi_1(X)&&\longrightarrow &&\pi_1(X_1)\cong Q^{+,-} \\
\endmatrix 
$$
of groups (see \cite{3} and \cite{5}). Iteration of this construction gives the second statement of the proposition. \qed
 \enddemo

Every diagram from Proposition 1  induces diagram (1.7). The surgery obstruction groups fitting into  the  diagram (1.7)  are well known \cite{10}:
$$
\matrix 
   i= & 0 & 1 & 2 & 3\\
   \\
 L_i(\Bbb Z/4^+)\cong & \Bbb Z^3 &0& \Bbb Z\oplus \Bbb Z/2 & \Bbb Z/2\\
   \\
    L_i(\Bbb Z/4^-)\cong & 0 &0 & \Bbb Z/2 & (\Bbb Z/2)^2\\
    \\
     L_i(Q)\cong & (\Bbb Z)^5 & |4| &  \Bbb Z/2 &  (\Bbb Z/2)^2 \\
     \\    
     L_i(Q^{+,-})\cong & |4| &0 & \Bbb Z\oplus \Bbb Z/2 & |4|, \\
     \\
     \endmatrix
   $$
where $|4|$ denotes a two-group of order 4. 

Let $n=4k+1$ and $\pi_1(X^n)\cong Q$. In this case commutative diagram (2.3) induces a diagram (see \cite{1})
$$
\matrix 0=L_1(\Bbb Z/4^+)& \overset{0}\to{\to} &   L_{1}(Q^+) &\cr
                &                     &\nearrow \ \ \ \ \ \ \searrow^p&\cr
                 &        LP_{0}(F_0)              &\Gamma\downarrow   \ \  \     &  L_{1}(\Bbb Z/4^+\to Q^+) \cr
                  &                     &\searrow \ \ \ \ \  \ \  
                  \nearrow_{r} &   \cr
0=LN_0(\Bbb Z/4^+\to Q^+)& \overset{0}\to{\to}  &  L_{0}(
Q^{+,-})&\cr
                       &                     &     ||           &\cr
0=L_0(\Bbb Z/4^-)& \overset{0}\to{\to} &   L_{0}(Q^{+,-}) &\cr
                &                     &\nearrow \ \ \ \ \ \searrow^{p}&\cr
                 &              LP_{3}(F_1)         &\Gamma\downarrow     \ \ \                &  L_{0}(\Bbb Z/4^-\to Q^{+,-}) \cr
                  &                     &\searrow \ \ \ \ \   \nearrow_{r} &   \cr
0=LN_3(\Bbb Z/4^-\to Q^{+,-})& \overset{0}\to{\to}  &  L_{3}(
Q^{-,-})&\cr
                       &                     &     ||           &\cr
\Bbb Z/2\cong L_3(\Bbb Z/4^+)& \overset{\phi (mono)}\to{\longrightarrow }
& L_{3}(Q^{-,-}) &\cr
                &                     &\nearrow \ \ \ \ \ \searrow&\cr
                 &              LP_{2}(F_2)         &\Gamma\downarrow  \ \ \                   &  L_{3}(\Bbb Z/4^+\to Q^{-,-}) \cr
                  &                     &\searrow \ \ \ \ \   \nearrow  &   \cr
LN_2(\Bbb Z/4^+\to Q^{-,-})& \overset{}\to{\to}  &  L_{2}(
Q^{+,+}). &\cr
\endmatrix
\tag 2.8
 $$
 where the groups  $LN_0(\Bbb Z/4^+\to Q^+)=0$ and  $LN_3(\Bbb Z/4^-\to Q^{+,-})=0$ were computed in \cite{1}. 
  
 All  maps $p$ and $r$  in  (2.8) are   monomorphisms, and   the surgery obstruction $\sigma(f)\in L_1(Q^+)$ 
 has the following properties   \cite{1}:
$$
p(\sigma(f))=r(\sigma(g)) \ \  \text{for some  element} \ \ \sigma(g)\in  L_{0}(Q^{+,-}),
$$
$$
p(\sigma(g))=r(\sigma(h)) \ \  \text{for some element} \ \ \sigma(h)\in  L_{3}(Q^{-,-}),
$$
and
$$
 \sigma(h) =\phi(a) 
$$
where  $a\in  L_{3}(\Bbb Z/4^{+})\cong \Bbb Z/2$ is  the nontrivial element. 
\bigskip

Note, that the results of \cite{1} immediately imply that 
$\Gamma^{3}(\sigma(f))=0$ for the diagram (2.3).

\smallskip

\subhead 3. Computing $\Gamma^{i}$ for Browder-Livesay filtrations
\endsubhead
\smallskip

In this section we prove the following result: let $X^{4k+1} \ (4k+1\geq 9)$ be  a manifold with $\pi_1(X)\cong Q$.  Consider the
Browder-Livesay filtration of $X$  with  a diagram of inclusion (2.3) -- (2.7). Then 
$\Gamma^3(\sigma(f))=0$ for the Cappell-Shaneson surgery obstruction  $\sigma(f)\in
L_1(Q^+)$. 

 The statement of the theorem for diagram (2.3) follows 
from diagram (2.8) and results of \cite{1}.  We shall give  
the statement of the theorem for other diagrams (2.4) -- (2.7) in 
 Theorems  3 and  4  of this section.
 
We introduce the following notations (see \cite{10}):
$$
\matrix
T^{-,-}=\Bbb Z[Q^{-,-}], \ \ \hat{T}_2^{-,-}=\hat{\Bbb Z}_2[Q^{-,-}],\\
R^{+}=\Bbb Z[\Bbb Z/4^{+}], \ \ \hat{R}_2^{+}=\hat{\Bbb Z}_2[\Bbb Z/4^+],
\endmatrix
$$
and similarly for other orientations.
Denote by 
$
L^{rel}_n(T^{-,-})$ the relative groups $ L_n(T^{-,-}\to
{T}_2^{-,-})$ fitting into the  relative exact sequence
$$
\to L_n(T^{-,-}) \to L_n(\hat{T}_2^{-,-}) \to L_n(T^{-,-}\to
{T}_2^{-,-}) \to L_{n-1}(T^{-,-}) \to 
$$
and similarly for other group rings and orientations.
>From the isomorphisms $Q^{+,-}\cong Q^{-,+}\cong Q^{-,-}$ which were described in Section 2, we obtain isomorphisms
$$
\matrix
T^{-,-}\cong T^{-,+}\cong T^{+,-},\\
\hat{T}_2^{-,-} \cong \hat{T}_2^{+,-} \cong \hat{T}_2^{-,+}\\
\endmatrix
$$
and isomorphisms 
$$
L_n^Y (T^{-,-}\to
{T}_2^{-,-})\cong L_n^Y (T^{+,-}\to
{T}_2^{+,-})\cong L_n^Y (T^{-,+}\to
{T}_2^{-,+}).
$$ 
We shall consider the
$L$-groups with decorations "prime". For the   case  of a group ring
$\Bbb Z[\pi]$ these $L$-groups coincide with surgery obstruction groups. 
We have isomorphisms (see \cite{1} and \cite{10}) 
$$
LN_n(\Bbb Z/4^+\to Q^{+,+})\cong L_n\left(\Bbb Z[\Bbb Z/4],\operatorname{Id}, -t^2\right)=L_n(B)
$$
and  the relative exact sequence
$$
\to L_n\left(\Bbb Z[\Bbb Z/4],\operatorname{Id}, -t^2\right)\to  L_n(\hat{\Bbb
Z}_2[\Bbb Z/4],\operatorname{Id}, -t^2)\to L^{rel}_0\left(\Bbb Z[\Bbb Z/4],\operatorname{Id},
-t^2\right)\to 
$$
Denote 
$$
L_n(B)=L_n\left(\Bbb Z[\Bbb Z/4],\operatorname{Id}, -t^2\right)
$$
and 
$$
L_n(\hat{B}_2)=
L_n\left(\hat{\Bbb Z}_2[\Bbb Z/4],\operatorname{Id}, - t^2\right).
$$
By \cite{1},  we have 
$$
LN_n(\Bbb Z/4^-\overset{j}\to{\to} Q^{-,-})\cong L_n\left(\Bbb Z[\Bbb Z/4], D, \pm t^2\right) =
L_n(A)
$$
where $D(t)=-t$. 
   In particular, we have an isomorphism 
$$
LN_n(\Bbb Z/4^-\overset{j}\to{\to} Q^{-,-})\cong LN_{n+2}(\Bbb Z/4^-\overset{j}\to{\to} Q^{-,-}).
$$
Denote 
$$
L_n(A)=L_n\left(\Bbb Z[\Bbb Z/4], D, \pm t^2\right) =
$$
and 
$$
L_n(\hat{A}_2)=
L_n\left(\hat{\Bbb Z}_2[\Bbb Z/4], D, \pm t^2\right).
$$
These groups fit into the  relative exact sequence
$$
\to L_n\left(\Bbb Z[\Bbb Z/4],D, \pm t^2\right)\to  L_n(\hat{\Bbb
Z}_2[\Bbb Z/4], D, \pm t^2)\to L^{rel}_n\left(\Bbb Z[\Bbb Z/4],D,
\pm t^2\right)\to 
$$
Braid of exact sequences  (1.5) also exists for relative groups (see \cite{6} and \cite{10}).  For
the inclusions $\Bbb Z/4^{\pm}\to Q^{+,{\pm}}$,   all relative groups and all  groups for the the  group ring over the ring $\hat{\Bbb Z}_2$  
are known \cite{10}.  In particular, we have 
$$
\matrix
n & & = & 0 & 1 & 2 & 3 \\
\\
L_n^{rel}(R^+)=&L_n^Y (R^{+}\to
{R}_2^{+})& \cong & 0 & \Bbb Z^3\oplus (\Bbb Z/2)^2 & 0 & \Bbb Z\\
\\
L_n^{rel}(T^+)=&
L_n^Y (T^{+}\to
{T}_2^{+}) &\cong & 0 & \Bbb Z^5\oplus (\Bbb Z/2)^4 & \Bbb Z/2 & (\Bbb Z/2)^2\\
\\
L_n^{rel}(B)=&
L_n^Y(B\to \hat{B}_2)&\cong& \Bbb Z/2& (\Bbb Z/2)^2& 0& (\Bbb Z)^2\oplus (\Bbb Z/2)^2\\
\\
L_n^{rel}(R^-)=&L_n^Y (R^{-}\to
{R}_2^{-})& \cong & \Bbb Z/2 & \Bbb Z/2\oplus \Bbb Z/2 & 0 & 0\\
\\
L_n^{rel}(T^-)=&
L_n^Y (T^{-,-}\to
{T}_2^{-,-}) &\cong & \Bbb Z/2 & \Bbb Z/2\oplus \Bbb Z/2 & 0 & \Bbb Z\\
\\
L_n^{rel}(A)=&
L_n^Y(A\to \hat{A}_2)&\cong& 0& \Bbb Z& 0& \Bbb Z\\
\endmatrix
\tag 3.1
$$

Consider the diagram 
$$
\matrix 0=L_1(\Bbb Z/4^+)& \overset{0}\to{\to} &   L_{1}(Q^+) &\cr
                &                     &s\nearrow \ \ \ \ \ \searrow^{\cong}&\cr
                 &                     &\Gamma\downarrow                    &  L_{1}(\Bbb Z/4^+\to Q^+) \cr
                  &                     &q\searrow \ \ \ \ \   ^{\cong}\nearrow  &   \cr
0=LN_0(\Bbb Z/4^+\to Q^+)& \overset{0}\to{\to}  &  L_{0}(
Q^{+,-})&\cr
                       &                     &     ||           &\cr
L_0(\Bbb Z/4^+)& \overset{Im = \Bbb Z/2}\to{\to} &   L_{0}(Q^{+,-})
&\cr
                &                     &s\nearrow \ \ \ \ \ \searrow^{}&\cr
                 &                     &\Gamma\downarrow                    &  L_{0}(\Bbb Z/4^+\to Q^{+,-})=\Bbb Z/2 \cr
                  &                     &q\searrow \ \ \ \ \   ^{}\nearrow  &   \cr
\Bbb Z/2=LN_3(\Bbb Z/4^+\to Q^{+,-})& \overset{mono}\to{\to}  &
L_{3}( Q^{+,+})=(\Bbb Z/2)^2&\cr
                       &                     &     ||           &\cr
L_3(\Bbb Z/4^+)& \overset{mono}\to{\to} &   L_{3}(Q^{+,+}) &\cr
                &                     &s\nearrow \ \ \ \ \ \searrow^{epi}&\cr
                 &                     &\Gamma\downarrow                    &  L_{3}(\Bbb Z/4^+\to Q^{-,-})=\Bbb Z/2 \cr
                  &                     &q\searrow \ \ \ \ \   {}\nearrow  &   \cr
LN_2(\Bbb Z/4^+\to Q^{+,+})& \overset{}\to{\to}  &  L_{2}(
Q^{+,-}).&\cr
\endmatrix
\tag 3.2 
 $$
 which is defined by diagram (2.4). Diagram (3.2) is similar to diagram (2.8). 
In the group  $L_{0}(Q^{+,-})$ lies the element $\sigma(g)$ which is
the image of the element $\sigma(f)\in L_{1}(Q^+)$ under the upper map
$\Gamma$ in (3.1) as follows from \cite{1}.

\proclaim{Theorem  3} 
 In  diagram (3.2), the element
 $\sigma(g)\in L_0(Q^{+,-})$ lies  in the image of the map
$$
 L_0(\Bbb Z/4^+) \longrightarrow
 L_{0}(Q^{+,-}),
 $$
and hence  $\Gamma^2(\sigma(f))=0$. From this  it  immediately follows that   $\Gamma^2(\sigma(f))=0$ for diagrams  (2.4) and (2.5). 
\endproclaim
\demo{Proof} Follows by  Lemma 1 and Lemma 2.

\proclaim{Lemma 1} 
Let $a$ be the nontrivial element of the relative group $L_0(T^{-,-}
\to \hat{T}_2^{-,-})=\Bbb Z/2$ and
$$
\partial \colon L_0(T^{-,-}
\to \hat{T}_2^{-,-})\to L_3(T^{-,-})
$$
be a map from the corresponding  relative exact sequence.
 Then
$$
\partial (a)
=\sigma(h)=\phi(x), \ \ \text{where} \ \ \ x\in L_3(\Bbb Z/4^+), \ \ x\ne 0.
$$
\endproclaim
\demo{Proof} 
Consider the commutative diagram 
$$
\CD \Bbb Z/4^+ @>>>  Q^{+,+}\\
@VVV  @VVV \\
\Bbb Z/2 @>>> \Bbb Z/2^+\oplus \Bbb Z/2^+
\endCD
\tag 3.3
$$
in which vertical maps are natural projections. By 
\cite{10}, the vertical maps in (3.3) induce  a retraction of
the braid of exact sequences  of relative groups for the  top
inclusion in  (3.3) to the corresponding braid of exact sequences  of
relative groups for the  bottom  inclusion in (3.3). 
 
As the kernel of
this retraction we obtain  a braid  of exact sequences with known groups and easily computed maps. 
All relative groups for the inclusion
$
\Bbb Z/2^+ \to \Bbb Z/2^+\oplus \Bbb Z/2^+
$ 
are known from \cite{10}:
$$
\matrix
n &  = & 0 & 1 & 2 & 3 \\
\\
L_n(\Bbb Z[\Bbb Z/2]\to \hat{\Bbb Z}_2[\Bbb Z/2] & \cong & 0 & (\Bbb Z\oplus \Bbb Z/2)^2 & 0 & 0\\
\\
L_n(\Bbb Z[(\Bbb Z/2)^2]\to \hat{\Bbb Z}_2[(\Bbb Z/2)^2] & \cong & 0 & (\Bbb Z\oplus \Bbb Z/2)^4& 0 & 0.\\
\endmatrix
$$
 All relative groups  for
the inclusion $\Bbb Z/4^+\to Q^{+,+}$
were given in  (3.1) (see \cite{10}).   From this we conclude that 
the map
$$
\Bbb Z/2=  L^{rel}_0\left(\Bbb Z[\Bbb Z/4],\operatorname{Id},
-t^2\right)\overset{\xi}\to{\longrightarrow} L^{rel}_0(T^{+,-})=\Bbb
Z/2 
\tag 3.4
$$
of relative groups is an isomorphism.  
The map $\xi$ from (3.4) lies in the commutative diagram 
$$
\CD
 0=L_0\left(\Bbb Z[\Bbb Z/4],\operatorname{Id}, -t^2\right)@>>>   L_0(\hat{\Bbb Z}_2[\Bbb Z/4],\operatorname{Id},
-t^2)@>{mono}>> L^{rel}_0\left(\Bbb Z[\Bbb Z/4],\operatorname{Id}, -t^2\right)=\Bbb Z/2 \\
@V0VV @VVV @V{\cong}V{\xi}V\\
|4|\cong L_0(T^{+,-})@>{epi}>>   L_0(\hat{T}_2^{+,-})(=\Bbb Z/2)@>0>> L^{rel}_0(T^{+,-})=\Bbb Z/2 \\
\endCD
\tag 3.5
$$
in which top row is a part of relative exact sequence for an inclusion 
$B\to \hat{B}_2$,  the
left group in top row is trivial by \cite{1}, the
bottom row is the part of the relative exact  sequence for an inclusion  $T^{+,-}
\to \hat{T}_2^{+,-}$, and the group and the maps in the bottom row were described
in \cite{10}. Since the map $\xi$ is an isomorphism, we obtain from commutativity 
of diagram (3.5) that 
 $$
L_0(\hat{\Bbb Z}_2[\Bbb Z/4],\operatorname{Id}, -t^2)=0.
$$
Consider a part 
$$
\CD \Bbb Z/2=L_3(\hat{R}_2^+) @>>> L_3(\hat{T}_2^{-,-})(= \Bbb Z/2)
@>>>L_3(\hat{\Bbb Z}_2[\Bbb Z/4], \operatorname{Id}, -t^2)@>>> L_1(\hat{T}_2^{+,+})\\
 @V{\Gamma}VV @V{\Gamma}VV @V{\Gamma}VV @V{\Gamma}VV\\
0=L_0(\hat{\Bbb Z}_2[\Bbb Z/4],\operatorname{Id}, -t^2) @>>>
L_2(\hat{T}_2^{+,+})(= \Bbb Z/2)
@>>>L_2(\hat{R}_2)@>{\cong}>> L_2(\hat{T}_2^{-,-})\\
\endCD
\tag 3.6
$$
of the two-row diagram for the inclusion
$
\hat{R}_2^+\to \hat{T}_2^{-,-}
$
which corresponds to the inclusion $\Bbb Z/4^+\to Q^{-,-}$ given by the map 
$t\to xy$.
The  rows of diagram (3.6) are chain complexes
with isomorphic homology groups. The right bottom horizontal map is
an isomorphism $\Bbb Z/2 \to \Bbb Z/2 $ (Arf-invariant), hence the
homology group in the member  $L_2(\hat{T}_2^{+,+})$ is $\Bbb Z/2$. Hence in top row the homology group in the member 
$L_3(\hat{T}_2^{-,-})$ is $\Bbb Z/2$  and 
hence 
the induced by the
inclusion  $\Bbb Z/4^+\to Q^{-,-} \ (t\to xy) $ map
$$
\alpha: \Bbb Z/2=L_3(\hat{R}_2^{+})\to L_3(\hat{T}_2^{-,-})=\Bbb Z/2
\tag 3.7
$$
is trivial.

The inclusion $\Bbb Z/4^+\to Q^{-,-} \ (t\to xy) $
induces the map of relative exact sequences, and we obtain the
commutative diagram
$$
\CD
 0= L_0(R^{+} \to \hat{R}_2^{+})@>0>>
L_3(R^{+})=\Bbb Z/2 @>{\cong}>> L_3(\hat{R}_2^{+})=\Bbb Z/2 \\
 @VVV @V{\phi}V{mono}V @V0VV\\
\Bbb Z/2= L_0(T^{-,-} \to \hat{T}_2^{-,-})@>{mono}>{\partial}>
 L_3(T^{-,-})=|4|@>{epi}>>
 L_3(\hat{T}_2^{-,-})=\Bbb Z/2\\
 \endCD
 \tag 3.8
$$
where the right hand  vertical map is trivial as we have proved. Now the statement of Lemma 1 follows from diagram (3.8). 
\qed
\enddemo

For the inclusion $\Bbb Z/4^-\to Q^{+,-} \ (t\to y) $ we consider
the following part of diagram (3.2)
$$
\matrix
  L_{0}(Q^{+,-}) && \cr
               &\searrow^{p}&\cr
                &&   L_{0}(\Bbb Z/4^-\to Q^{+,-}) \cr
                   & ^{r}\nearrow  &   \cr
 L_{3}(Q^{-,-})&&\cr
\endmatrix
\tag 3.9
 $$
 where  $\sigma(g)\in L_{0}(Q^{+,-})$ and $\sigma(h)\in
 L_{3}(Q^{-,-})$. We have a natural map of the corresponding diagram
 of relative groups to diagram (3.2) and hence a map of a diagram of relative groups to diagram  (3.9). Thus we obtain the commutative diagram 
$$
\CD
 0\ne b\in \Bbb Z/2\oplus \Bbb Z/2 =L_1(T^{+,-} \to \hat{T}_2^{+,-}) @>{Im =\Bbb Z/2}>{\delta}>
 L_{0}(Q^{+,-})\backepsilon \sigma(g)
 \\
  @V{epi}VV   @V{p}V{mono}V\\
 \Bbb Z/2= L_1\left(\matrix R^-&\to& \hat{R}_2^{-}\\
                   \downarrow && \downarrow\\
                   T^{+,-}& \to&  \hat{T}_2^{+,-}\\
                   \endmatrix\right) @>{mono}>>    L_{0}(\Bbb Z/4^-\to Q^{+,-})
                   \\
@A{\cong}AA @A{r}A{mono}A\\
 0\ne a\in \Bbb Z/2 =L_0(T^{-,-} \to \hat{T}_2^{-,-})@>{mono}>{\partial}>
  L_{3}(Q^{-,-}) \backepsilon \sigma(h)\\
\endCD
\tag 3.10
$$
in which the maps $p$ and $r$  are monomorphisms by  \cite{1}.
The middle horizontal map in (3.10) is a monomorphism since the diagonal
map
$$
L_0(T^{-,-} \to \hat{T}_2^{-,-})\to L_{0}(\Bbb Z/4^-\to Q^{+,-})
$$
from bottom square is a monomorphism  as follows from Lemma 1. The left hand upper vertical map in (3.10) is an epimorphism  by \cite{10}. Hence from commutativity of (3.10) we obtain 
 that there exists an element
 $$
 b\in L_1(T^{+,-} \to \hat{T}_2^{+,-})
 $$
 such that $\delta(b)= \sigma(g)\in  L_{0}(Q^{+,-})$.

Consider the map of relative exact sequence of  the inclusion  $R^+\to \hat{R}_2^+$
to relative exact sequence of the inclusion  $T^{+,-}\to \hat{T}_2^{+,-}$ which is 
induced by the inclusion $\Bbb Z/4^+\to Q^{+,-}$. We can write down
the commutative diagram 
$$
\CD
 \Bbb Z^3\oplus (\Bbb Z/2)^3=L_1(R^+\to \hat{R}_2^{+}) @>{epi}>>
 L_0(R^+)  =\Bbb Z^3 @.\\
@V{Im=\Bbb Z/2}VV @V{i_*}VV @.\\
0\ne b\in (\Bbb Z/2)^2=L_1(T^{+,-}\to \hat{T}_2^{+,-}) @>{Im =\Bbb
Z/2}>{\delta}>
 L_0(T^{+,-})=|4|@. (\backepsilon \sigma(g)=\delta(b)) \\
\endCD
\tag 3.11
$$
where the left hand vertical map follows from \cite{10}.

\proclaim{Lemma 2} The map
$$
\Bbb Z^3=L_0(R^+)\overset{{i_*}}\to{\to} L_0(T^{+,-})
$$
in diagram (3.11) has the image $\Bbb Z/2$ and hence the element
$\sigma(g)$ lies in the image of this map.
\endproclaim
\demo{Proof} Commutative diagram (3.3) induces a commutative diagram 
$$
\CD
LN_1(\Bbb Z/4^+\to  Q^{+,+})@>>> L_3(Q^{+,+})=(\Bbb Z/2)^2\\
@V{epi}VV  @V{mono}VV \\
LN_1(\Bbb Z/2^+\to  \Bbb Z/2^+\oplus\Bbb Z/2^+) @>{mono}>> L_3(\Bbb
Z/2^+\oplus \Bbb Z/2^+)=(\Bbb Z/2)^3\\
@\vert @. \\
L_3(\Bbb Z/2^+)=\Bbb Z/2 @.
\endCD
\tag 3.12
$$
in which the right   vertical map is a monomorphism by \cite{1} and the left  vertical map is an epimorphism as follows
from  the long exact sequence of the pull-back diagram (see
\cite{1})
$$
\CD
(\Bbb Z[\Bbb Z/4], Id, -t^2)@>>> (\Bbb Z[i], Id, 1)\\
@VVV  @VVV \\
(\Bbb Z[\Bbb Z/2], Id, -1) @>>> (\Bbb Z/2[\Bbb Z/2], Id, 1),
\endCD
$$
where $LN_i(\Bbb Z/4^+\to  Q^{+,+})= L_i(\Bbb Z[\Bbb Z/4], Id,
-t^2)$. Thus  the image of the upper horizontal map in (3.12) is $\Bbb
Z/2$. Now consider the following part of the two-row diagram for the
inclusion $\Bbb Z/4^+\to  Q^{+,-}$
$$
\CD
0=L_1(Q^{+,-})@>>> LN_1(\Bbb Z/4^+\to  Q^{+,+})@>{Im=\Bbb Z/2}>> L_3(Q^{+,+})=(\Bbb Z/2)^2@>>> L_3(\Bbb Z/4^+)\\
@V{\Gamma}VV  @V{\Gamma}VV @V{\Gamma}VV  @V{\Gamma}VV \\
\Bbb Z^5=L_0(Q^{+,+})@>{coker=\Bbb Z/2}>> L_0(\Bbb Z/4^+)=\Bbb
Z^3@>{i_*}>>
L_0(Q^{+,-})=|4|@>>> LN_0=0\\
\endCD
\tag 3.13
$$
in which rows are chain complexes with isomorphic homology groups
and vertical maps $\Gamma$ are isomorphisms of homologies. The map
$$
\Bbb Z^5=L_0(Q^{+,+})\to L_0(\Bbb Z/4^+)=\Bbb Z^3
$$
has cokernel $\Bbb Z/2$ as follows from consideration of the
corresponding diagram of relative groups. It follows now that image
${i_*}$ in (3.13) is $\Bbb Z/2$ and the lemma is proved.
 \qed \enddemo
 Now Theorem 3 is proved.
 \enddemo

On the base of diagram (2.6) we can write down the following  commutative diagram 
of surgery obstruction groups

$$
\matrix0=L_1(\Bbb Z/4^+)& \overset{0}\to{\to} &   L_{1}(Q^+) &\cr
                &                     &\nearrow \ \ \ \ \ \ \searrow {p}&\cr
                 &                     &\Gamma\downarrow   \ \  \     &  
L_{1}(\Bbb Z/4^+\to Q^+) \cr
                  &                     &\searrow \ \ \ \ \  \ \  \nearrow {r} &   \cr
0=LN_0(\Bbb Z/4^+\to Q^+)& \overset{0}\to{\to}  &  L_{0}(Q^{+,-})&\cr
                       &                     &     ||           &\cr
0=L_0(\Bbb Z/4^+)& \overset{j_*=0}\to{\longrightarrow} &   L_{0}(Q^{+,-}) &\cr
                &                     &\nearrow \ \ \ \ \ \searrow{p}&\cr
                 &                     &\Gamma\downarrow                    & 
 L_{0}(\Bbb Z/4^-\to Q^{+,-}) \cr
                  &                     &\searrow \ \ \ \ \   \nearrow{r}  &   \cr
LN_3(\Bbb Z/4^-\to Q^{+,-})& \overset{0}\to{\to}  &  L_{3}( Q^{-,-})&\cr
                       &                     &     ||           &\cr
\Bbb Z/2\oplus \Bbb Z/2\cong L_3(\Bbb Z/4^-)& \overset{j_*}\to{\longrightarrow }
& L_{3}(Q^{-,-}) &\cr
                &                     &\nearrow \ \ \ \ \ \searrow&\cr
                 &                     &\Gamma\downarrow                    & 
 L_{3}(\Bbb Z/4^-\to Q^{-,-})\cr
                  &                     &\searrow \ \ \ \ \   \nearrow  &   \cr
LN_2(\Bbb Z/4^-\to Q^{-,-})& \overset{}\to{\longrightarrow}  &  L_{2}(
Q^{+,-})&\cr
\endmatrix
\tag 3.14
 $$
which is similar to diagram (2.8). By \cite{1}, 
 $\sigma(g)\in L_{0}(Q^{+,-})$,   $\sigma(h)\in L_{3}(Q^{-,-})$, and 
 $p(\sigma(g))=r(\sigma(h))$.

\proclaim{Theorem 4} The map 
$$
\Bbb Z/2\oplus \Bbb Z/2=L_3(\Bbb Z/4^-) \overset{j_*}\to{\longrightarrow }
 L_{3}(Q^{-,-})
$$
in  diagram (3.14) is an isomorphism. 
Hence   $\Gamma^3(\sigma(f))\in L_2(Q^{+,-})$ is trivial for diagram (3.14). From this it  follows 
   that  $\Gamma^3(\sigma(f))=0$ for diagrams (2.6) and (2.7).
\endproclaim
\demo{Proof}
For the inclusion 
$j:\Bbb Z/4^-\to Q^{+,-}$, consider a two-row  diagram of relative $L$-groups, 
which is part of the braid of exact sequences of relative $L$-groups 
$$
\CD
 \Bbb Z/2 @. \Bbb Z/2 @. 0 @. 0 @. 0 @. 0\\
  @\vert @\vert @\vert @\vert @\vert @\vert\\
\overset{j^!}\to{\to} L_0^{rel}(R^{-})@>>>L_0^{rel}(T^{-}) @>>> L_2^{rel}(A)@>>> L_2^{rel}(T^{-}) @>>> L_2^{rel}(R^{-})@>>>L_2^{rel}(T^-)\to\\
 @V{\Gamma}VV  @V{\Gamma}VV @V{\Gamma}VV  @V{\Gamma}VV @V{\Gamma}VV @V{\Gamma}VV\\
\to L_3^{rel}(A)@>>>L_3^{rel}(T^{-}) @>>>
L_3^{rel}(R^{-})@>>>L_3^{rel}(T^{-}) @>>>L_1^{rel}(A) @>>> L_1^{rel}(T^-)\to\\
   @\vert @\vert @\vert @\vert @\vert @\vert\\
 \Bbb Z @. \Bbb Z @. 0 @. \Bbb Z @. \Bbb Z @. (\Bbb Z/2)^2\\
\endCD
\tag 3.15
$$
$$
\CD
 @. 0 @. \Bbb Z/2 \\
 @. @\vert @\vert \\
@>>> L_0^{rel}(A)@>>>L_0^{rel}(T^{-})\overset{j^!}\to{\to}\\
@. @V{\Gamma}VV  @V{\Gamma}VV \\
@>>> L_1^{rel}(R^-)@>{j_*}>>L_1^{rel}(T^{-})\to\\
@.   @\vert @\vert \\
@. (\Bbb Z/2)^2 @. (\Bbb Z/2)^2\\
\endCD
$$
 Now we shall prove Lemmas 3,4,5, and 6 which yields  a proof of Theorem 
 4. 
 \proclaim{Lemma 3}  In diagram (3.15),  the map
 $$
\Bbb Z=L_3^{rel}(A)@>>>L_3^{rel}(T^{-})=\Bbb Z
$$
is isomorphism, 
and the map 
$$
\Bbb Z=L_3^{rel}(T^{-})\to L_1^{rel}(A)=\Bbb Z
$$
is a multiplication by 2 with coimage $\Bbb Z/2$.
\endproclaim
\demo{Proof}  Suppose that the map 
$$
\Bbb Z/2=L_0^{rel}(T^{-})\overset{j^!}\to{\longrightarrow} L_0^{rel}(R^-)=\Bbb Z/2
$$
in diagram (3.15) is an isomorphism. Then, comparing homology of top and bottom rows of diagram (3.15), 
we obtain a contradiction. Hence the map 
$$
\CD
\Bbb Z/2=L_0^{rel}(R^-)@>{j_*}>>L_0^{rel}(T^{-})=\Bbb Z/2
\endCD
$$
in top row of (3.15) is an isomorphism. From this the result follows by 
diagram chasing in diagram  (3.15). 
\qed
\enddemo

For an  inclusion 
$
j: \Bbb Z/4^- \to  Q^{-}, 
$
we have  a  two-row diagram of $L$-groups
$$
\CD
  0 @. |4| @. ? @. \Bbb Z\oplus\Bbb Z/2  @.\Bbb Z/2 @.\Bbb Z\oplus\Bbb Z/2\\
  @\vert @\vert @\vert @\vert @\vert @\vert\\
\overset{}\to{\to} L_0(\Bbb Z/4^{-})@>>>L_0(Q^{-}) @>>> LN_2@>>> L_2(Q^{-}) @>{j^!}>> L_2(\Bbb Z/4^{-})@>{j_*}>>L_2(Q^-)\to\\
 @V{\Gamma}VV  @V{\Gamma}VV @V{\Gamma}VV  @V{\Gamma}VV @V{\Gamma}VV @V{\Gamma}VV\\
\to LN_3@>>>L_3(Q^{-}) @>>>
L_3(\Bbb Z/4^{-})@>{j_*}>>L_3(Q^{-}) @>>>LN_1 @>>> L_1(Q^-)\to\\
   @\vert @\vert @\vert @\vert @\vert @\vert\\
0 @. |4| @. \Bbb Z/2\oplus\Bbb Z/2 @. |4| @. 0 @. 0\\
\endCD
\tag 3.16
$$
$$
\CD
 @. ? @. |4| \\
 @. @\vert @\vert \\
@>>> LN_0 @>>>L_0(Q^{-})\overset{j^!}\to{\to}\\
@. @V{\Gamma}VV  @V{\Gamma}VV \\
@>>> L_1(\Bbb Z/4^-)@>{}>>L_1(Q^{-})\to\\
@.   @\vert @\vert \\
@. 0 @. 0\\
\endCD
$$

 \proclaim{Lemma 4}  In diagram (3.16),  the map
 $$
\Bbb Z/2=L_2(\Bbb Z/4^-)@>{j_*}>>L_2(Q^{-})=\Bbb Z\oplus\Bbb Z/2
$$
is a monomorphism, and hence the map 
$$
\CD 
L_2(Q^{-}) @>{j^!}>> L_2(\Bbb Z/4^{-})\\
\endCD
$$
is trivial. 
\endproclaim
\demo{Proof} Consider the commutative diagram of oriented groups
$$
\CD 
\Bbb Z/4^- @>j(t)=y>> Q^{+,-}\\
@V{t^2= 1}VV @V{x^2=y^2=1}VV\\
\Bbb Z/2^-@>>> \Bbb Z/2^+\oplus \Bbb Z/2^-.\\
\endCD
\tag 3.17
$$
Diagram (3.17) induces the following diagram of $L$-groups 
$$
\CD 
\Bbb Z/2=L_2(\Bbb Z/4^-) @>{j_*}>> L_2(Q^{+,-})\\
@V{\cong}VV @VVV\\
\Bbb Z/2=L_2(\Bbb Z/2^-)@>{\cong}>> L_2(\Bbb Z/2^+\oplus \Bbb Z/2^-)=\Bbb Z/2\\
\endCD
\tag 3.18
$$
in which the left vertical map is isomorphism (Arf-invariant is preserved), 
and the bottom horizontal map is an isomorphism  (the bottom horizontal map  is an inclusion of a direct summand). Now we obtain the monomorphism of top  row in (3.18). The result follows. \qed
\enddemo 

\proclaim{Lemma 5}  We have  isomorphisms \cite{6}
$$
LN_2(\Bbb Z/4^-\to Q_3^{-})\cong LN_0(\Bbb Z/4^-\to Q_3^{-})\cong \Bbb Z\oplus \Bbb Z/2
$$
and the map 
$$
\Bbb Z\cong L_{2n+1}^{rel}(A)\to L_{2n}(A)\cong LN_{2n}(\Bbb Z/4^-\to Q_3^{-})
$$
is an inclusion of a direct summand.
\endproclaim 
\demo{Proof} 
 Consider a commutative diagram 
$$
\CD
\Bbb Z=L_3^{rel}(A)@>{\cong}>> \Bbb Z=L_3^{rel}(T^-) \\
 @V{\partial}VV @V{\partial_1}VV  \\
  LN_2 @>>> \Bbb Z\oplus\Bbb Z/2=L_2(Q^-)\\
 \endCD
\tag 3.19
$$
in which vertical maps fit into the corresponding relative exact sequences. The right vertical map ${\partial_1}$ in (3.19)  is an inclusion of a direct summand
by \cite{10}, and upper horizontal map is an isomorphism by Lemma 3.
Hence the left vertical map in (3.19) is an inclusion of a direct summand. Now from  diagram 
(3.16) we obtain the diagram 

$$
\CD
   @. |4| @. ? @. \Bbb Z/2  @.\\
  @. @\vert @\vert @\vert @.\\
0@>>> L_0(Q^{-}) @>>> LN_2/\{{\operatorname Im}\ \partial\}@>>>L_2(Q^{-}) /\{{\operatorname Im}\ \partial_1\} @>>> 0\\
 @. @V{\Gamma}VV @V{\Gamma}VV  @V{\Gamma}VV @. \\
0 @>>>L_3(Q^{-}) @>>>
L_3(\Bbb Z/4^{-})@>{j_*}>>L_3(Q^{-}) @>>> 0\\
   @. @\vert @\vert @\vert @.\\
 @. |4| @. \Bbb Z/2\oplus\Bbb Z/2 @. |4| @.\\
\endCD
\tag 3.20
$$
in which the homology groups in the corresponding places of upper and bottom 
rows are isomorphic by isomorphisms $\Gamma$. Hence we have the isomorphism 
$$
LN_2/\{{\operatorname Im} \ \partial\}\cong \Bbb Z/2
$$
and 
$$
LN_2\cong \Bbb Z\oplus\Bbb Z/2.
$$
The isomorphism $LN_2\cong LN_0$ follows from 2-periodicity. \qed
\enddemo

\proclaim{Lemma 6} The map 
$$
\Bbb Z\oplus \Bbb Z/2=LN_2\to L_2(Q^{+,-})=\Bbb Z\oplus \Bbb Z/2
$$
in diagram (3.16) is an isomorphism
\endproclaim
\demo{Proof}
Consider the commutative diagram 
$$
\CD 
\Bbb Z\oplus \Bbb Z/2\cong LN_0(\Bbb Z/4^-\to Q^{-,-})@>{\cong}>> L_0\left(\Bbb Z[\Bbb Z/4], D, t^2\right)@>{epi}>> L_0(Q^{+,-})\\
@V{epi}VV @V{t^2=1}VV @V{x^2=y^2=1}VV\\
\Bbb Z/2\cong LN_0(\Bbb Z/2^-\to \Bbb Z/2^-\oplus \Bbb Z/2^-)@>{\cong}>> L_0(\Bbb Z/2^-)@>{\cong}>> L_0(\Bbb Z/2^+\oplus \Bbb Z/2^-)\\
\endCD
\tag 3.21
$$
in which the top right hand  horizontal  map is an epimorphism as follows from (3.16).   
Hence the left vertical map provides  an isomorphism of the torsion part $\Bbb Z/2$ of the group 
$LN_0$ on the group  $L_0(\Bbb Z/2^-)=\Bbb Z/2$ (the direct summand $\Bbb Z$ of the group
$LN_0$ maps trivially since it factors acroos the corresponding map of relative groups
which is trivial). 

Now we can write down a commutative  diagram  
$$
\CD 
\Bbb Z\oplus \Bbb Z/2\cong LN_2(\Bbb Z/4^-\to Q^{-,-})@>{\cong}>> L_2\left(\Bbb Z[\Bbb Z/4], D, t^2\right)@>>> L_2(Q^{+,-})\cong \Bbb Z\oplus\Bbb Z/2\\
@V{epi}VV @V{t^2=1}VV @V{x^2=y^2=1}VV\\
\Bbb Z/2\cong LN_2(\Bbb Z/2^-\to \Bbb Z/2^-\oplus \Bbb Z/2^-)@>{\cong}>> L_2(\Bbb Z/2^-)@>{\cong}>> L_2(\Bbb Z/2^+\oplus \Bbb Z/2^-)\\
\endCD
\tag 3.22
$$
in which 
the left  vertical map induces  an isomorphism of the torsion part $\Bbb Z/2$ of the group 
$LN_2$ on the group  $L_2(\Bbb Z/2^-)=\Bbb Z/2$ by  2-periodicity.
The direct summand $\Bbb Z$ of the group $LN_2$ maps isomorphically to the direct summand $\Bbb Z$ of the group $L_2(Q^{+,-})$ as follows from Lemma 3   and Lemma 5. The diagram (3.22) provides a similar result about torsion subgroups.  Lemma 6 is proved. Now   Theorem 4 follows.
\qed
\enddemo
\enddemo

\subhead  4. Braids of exact sequences
\endsubhead 
\bigskip  

In this section we give explicit results about braids of exact sequences of $L_*$-groups for various  index 2  inclusions $\rho \to Q$. The results of this section are direct corollaries of computations in Section 3 and
\cite{1}, \cite{6}, and \cite{10}.  

\proclaim{Theorem 5} 
For an  inclusion 
$
j: \Bbb Z/4^- \to  Q^{-}, 
$
we have  the following  two-row diagram of $L$-groups
$$
\CD
  0 @. (\Bbb Z/2)^2 @. \Bbb Z\oplus \Bbb Z/2 @. \Bbb Z\oplus\Bbb Z/2  @.\Bbb Z/2  \\
  @\vert @\vert @\vert @\vert @\vert \\
\overset{}\to{\to} L_0(\Bbb Z/4^{-})@>>>L_0(Q^{-}) @>{0}>> LN_2@>{\cong}>> L_2(Q^{-}) @>{j^!=0}>> L_2(\Bbb Z/4^{-})\to\\
 @V{\Gamma}VV  @V{\Gamma}VV @V{\Gamma}VV  @V{\Gamma}VV @V{\Gamma}VV \\
\to LN_3@>>>L_3(Q^{-}) @>{0}>>
L_3(\Bbb Z/4^{-})@>{\cong}>>L_3(Q^{-}) @>>>LN_1\to\\
   @\vert @\vert @\vert @\vert @\vert \\
0 @. (\Bbb Z/2)^2 @. \Bbb Z/2\oplus\Bbb Z/2 @. (\Bbb Z/2)^2 @. 0 \\
\endCD
\tag 4.1
$$
$$
\CD
@.\Bbb Z\oplus\Bbb Z/2 @. \Bbb Z\oplus \Bbb Z/2 @. (\Bbb Z/2)^2 \\
 @.@\vert @\vert @\vert \\
@>{j_*}>{mono}>L_2(Q^-)@>{Coker=(\Bbb Z/2)^2}>> LN_0 @>{epi}>>L_0(Q^{-})\overset{j^!=0}\to{\to}\\
@. @V{\Gamma}VV  @V{\Gamma}VV  @V{\Gamma}VV \\
@>>> L_1(Q^-)@>>> L_1(\Bbb Z/4^-)@>{}>>L_1(Q^{-})\to\\
@. @\vert  @\vert @\vert \\
@.0 @.  0 @. 0\\
\endCD
$$
\endproclaim 

\proclaim{Corollary 1} There exist  isomorphisms
$$
L_0(Q^-)\cong L_3(Q^-)\cong \Bbb Z/2\oplus \Bbb Z/2.
$$
\endproclaim 

\proclaim{Corollary 2} Let $Y\subset X$ be be a Browder-Livesay pair of manifolds with  a push-out  square of fundamental groups 
$$
\Psi=\CD 
\pi_1(\partial U)@>\cong >> \pi_1(X\setminus Y)\\
@VVV @VVV \\
\pi_1(Y)@>>> \pi_1(X)\\
\endCD\ 
= \
\CD 
\Bbb Z/4^-@>>> \Bbb Z/4^-\\
@VVV @VVV \\
Q^{-}@>>>Q^{-} \\
\endCD
$$
in the corresponding splitting problem.  Then we have the following  isomorphisms 
$$
LP_n(\Psi)= (\Bbb Z/2)^2,  \Bbb Z/2, \Bbb Z\oplus (\Bbb Z/2)^3, (\Bbb Z/2)^2
$$
for $n=0,1,2,3$ $\bmod 4$,  respectively. 
\endproclaim 

\proclaim{Theorem 6} 
For an  inclusion 
$
j: \Bbb Z/4^+ \to  Q^{+,-}, 
$
we have  the following  two-row diagram of $L$-groups
$$
\CD
  0 @. 0 @. \Bbb Z/2@. (\Bbb Z/2)^2  @.\Bbb Z/2 \\
  @\vert @\vert @\vert @\vert @\vert \\
\overset{}\to{\to} L_1(\Bbb Z/4)@>>>L_1(Q^-) @>{0}>> LN_1@>{mono}>> L_3(Q) @>{0}>> L_3(\Bbb Z/4)\to\\
 @V{\Gamma}VV  @V{\Gamma}VV @V{\Gamma}VV  @V{\Gamma}VV @V{\Gamma}VV \\
\to LN_2@>{mono}>>L_0(Q) @>{Coker=\Bbb Z/2}>>
L_0(\Bbb Z/4)@>{Im=\Bbb Z/2}>>L_0(Q^-) @>>>LN_0\to\\
   @\vert @\vert @\vert @\vert @\vert \\
\Bbb Z^2 @. \Bbb Z^5 @. \Bbb Z^3 @. (\Bbb Z/2)^2 @. 0 \\
\endCD
\tag 4.2
$$
$$
\CD
@. (\Bbb Z/2)^2 @.  (\Bbb Z/2)^2 @. (\Bbb Z/2)^2 \\
 @.@\vert @\vert @\vert \\
@>{mono}>>L_3(Q^-)@>{0}>> LN_3@>{\cong}>>L_1(Q)\overset{0}\to{\to}\\
@. @V{\Gamma}VV  @V{\Gamma}VV  @V{\Gamma}VV \\
@>>> L_2(Q)@>0>> L_2(\Bbb Z/4)@>{\cong}>>L_2(Q^-)\overset{0}\to{\to}\\
@. @\vert  @\vert @\vert \\
@. \Bbb Z/2 @.  \Bbb Z\oplus \Bbb Z/2@. \Bbb Z\oplus \Bbb Z/2\\
\endCD
$$
where 
$LN_i=LN_i(\Bbb Z/4\to Q^+)\cong LN_{i+2}(\Bbb Z/4\to Q^-)$.

\endproclaim 

\bigskip

\proclaim{Corollary 3} There exists  an isomorphism
$$
L_1(Q^+)\cong \Bbb Z/2\oplus \Bbb Z/2.
$$
\endproclaim 

\proclaim{Corollary 4} Let $Y\subset X$ be be a Browder-Livesay pair of manifolds with  a push-out  square of fundamental groups 
$$
\Psi_1=\CD 
\pi_1(\partial U)@>\cong >> \pi_1(X\setminus Y)\\
@VVV @VVV \\
\pi_1(Y)@>>> \pi_1(X)\\
\endCD\ 
= \
\CD 
\Bbb Z/4^+@>>> \Bbb Z/4^+\\
@VVV @VVV \\
Q^{+,+}@>>>Q^{+,-} \\
\endCD
$$
in the corresponding splitting problem.  Then we have the following  isomorphisms 
$$
LP_n(\Psi_1)= \Bbb Z^2, \Bbb Z\oplus (\Bbb Z/2)^3,  (\Bbb Z/2)^2, (\Bbb Z/2)^3
$$
for $n=0,1,2,3$ $\bmod 4$,  respectively. 
\endproclaim

\proclaim{Theorem 7} 
For an  inclusion 
$
j: \Bbb Z/4^+ \to  Q^{+}, 
$
we have  the following  two-row diagram of $L$-groups
$$
\CD
  0 @. (\Bbb Z/2)^2 @. (\Bbb Z/2)^2@. (\Bbb Z/2)^2  @.\Bbb Z/2 \\
  @\vert @\vert @\vert @\vert @\vert \\
\overset{}\to{\to} L_1(\Bbb Z/4)@>>>L_1(Q) @>{0}>> LN_3@>{\cong}>> L_3(Q^{-}) @>{0}>> L_3(\Bbb Z/4)\overset{mono}\to{\longrightarrow}\\
 @V{\Gamma}VV  @V{\Gamma}VV @V{\Gamma}VV  @V{\Gamma}VV @V{\Gamma}VV \\
\to LN_0@>>>L_0(Q^{-}) @>{0}>>
L_0(\Bbb Z/4)@>{mono}>>L_0(Q) @>{epi}>>LN_2\overset{0}\to{\to}\\
   @\vert @\vert @\vert @\vert @\vert \\
0 @. (\Bbb Z/2)^2 @. \Bbb Z^3 @. \Bbb Z^5 @. \Bbb Z^2 \\
\endCD
\tag 4.3
$$
$$
\CD
@. (\Bbb Z/2)^2 @.  \Bbb Z/2 @. 0 \\
 @.@\vert @\vert @\vert \\
@>{mono}>>L_3(Q)@>{0}>> LN_1@>{}>>L_1(Q^{-})\overset{}\to{\to}\\
@. @V{\Gamma}VV  @V{\Gamma}VV  @V{\Gamma}VV \\
@>0>> L_2(Q^-)@>{Coker=(\Bbb Z/2)^2}>{Ker=\Bbb Z/2}> L_2(\Bbb Z/4)@>{epi}>>L_2(Q)\to\\
@. @\vert  @\vert @\vert \\
@.\Bbb Z\oplus \Bbb Z/2 @.  \Bbb Z\oplus \Bbb Z/2@. \Bbb Z/2\\
\endCD
$$
where 
$LN_i=LN_i(\Bbb Z/4\to Q^+)\cong LN_{i+2}(\Bbb Z/4\to Q^-)$.
\endproclaim 

\bigskip

\proclaim{Corollary 5} Let $Y\subset X$ be be a Browder-Livesay pair of manifolds with  a push-out  square of fundamental groups 
$$
\Psi_2=\CD 
\pi_1(\partial U)@>\cong >> \pi_1(X\setminus Y)\\
@VVV @VVV \\
\pi_1(Y)@>>> \pi_1(X)\\
\endCD\ 
= \
\CD 
\Bbb Z/4^+@>>> \Bbb Z/4^+\\
@VVV @VVV \\
Q^{+,-}@>>>Q^{+,+} \\
\endCD
$$
in the corresponding splitting problem.  Then we have the following  isomorphisms 
$$
LP_n(\Psi_2)= (\Bbb Z/2)^2,  (\Bbb Z/2)^2,  (\Bbb Z/2)^2, \Bbb Z^3\oplus (\Bbb Z/2)^2
$$
for $n=0,1,2,3$ $\bmod 4$,  respectively. 
\endproclaim

\bigskip

Authors' addresses:

\bigskip

\noindent
Friedrich Hegenbarth: Dipartimento di Matematica, Universit\`a di
Milano,
Via Saldini n. 50, 20133 Milano, Italy
\smallskip
\noindent
E--mail: friedrich.hegenbarth\@unimi.it

\bigskip
\noindent
Yuri Muranov: Institute of Mathematics, Jagiellonian University,
ul. prof. Stanislawa Lojasiewicza 6, 30-348 Krakow, Poland
\smallskip
\noindent
E--mail: ymuranov\@mail.ru

\bigskip

\noindent Du\v san Repov\v s: Faculty of Mathematics ans  Physics,
and Faculty of Education,
University of
Ljubljana, Jadranska 19, 1000 Ljubljana, Slovenia;
\smallskip

\noindent E--mail: dusan.repovs\@guest.arnes.si

\enddocument
\bye